\documentclass{amsart}

\usepackage{color} % This package was added to write corrections in red

\usepackage{appendix}

\newtheorem{theorem}{Theorem}[section]
\newtheorem{prop}[theorem]{Proposition}
\newtheorem{proposition}[theorem]{Proposition}

\theoremstyle{definition}

\theoremstyle{remark}
\newtheorem{remark}[theorem]{Remark}

\newcommand{\bd}{{\partial}}
\newcommand{\bH}{\mathbb{H}}

\newcommand{\proj}{{\mathbb P}}

\newcommand{\cG}{{\mathcal{G}}}

\newcommand{\cN}{{\mathcal{N}}}

\newcommand{\cE}{{\mathcal{E}}}

\newcommand{\cR}{{\mathcal{R}}}

\newcommand{\R}{{\mathbb{R}}}
\newcommand{\Z}{{\mathbb{Z}}}

\newcommand{\N}{\mathbb{N}}

\renewcommand{\part}[2]{\frac{\partial #2}{\partial x_{#1}}}

\def\sideremark#1{\ifvmode\leavevmode\fi\vadjust{\vbox to0pt{\vss% the remark
 \hbox to 0pt{\hskip\hsize\hskip1em%                          will appear only
 \vbox{\hsize3cm\tiny\raggedright\pretolerance10000%          on the side
  \noindent #1\hfill}\hss}\vbox to8pt{\vfil}\vss}}}%
\begin{document}

\title[Conformal invariants from nullspaces]{Nullspaces of Conformally Invariant Operators. Applications to $Q_{k}$-curvature}

\author[Y. Canzani]{Yaiza Canzani}
\address{Department of Mathematics and Statistics, McGill University, Montr\'eal,
Ca\-na\-da.} 
\email{canzani@math.mcgill.ca}

\author[A.R. Gover]{A. Rod Gover}
\address{Department of Mathematics, University of Auckland, New Zealand \&
Mathematical Sciences Institute, Australian National University, Canberra, Australia.}
\email{gover@math.auckland.ac.nz}

\author[D. Jakobson]{Dmitry Jakobson}
\address{Department of Mathematics and
Statistics, McGill University, Montr\'eal, Ca\-na\-da.}
\email{jakobson@math.mcgill.ca}

\author[R. Ponge]{Rapha\"el Ponge}
\address{Department of Mathematical Sciences, Seoul National University, Seoul, Korea.} 
\email{ponge.snu@gmail.com}



\thanks{YC was supported by Schulich Fellowship of McGill University (Canada). 
ARG was supported by Marsden Grant 10-UOA-113 (New Zealand). DJ was
  supported by NSERC and FQRNT grants and Dawson Fellowship of McGill University (Canada).  RP was supported by a JSPS Grant-in-Aid (Japan) and Research Resettlement Fund of Seoul National University (Korea).}  

\begin{abstract}
We study conformal invariants that arise from functions in the
nullspace of conformally covariant differential operators.
The invariants include nodal sets and the topology of nodal domains
of eigenfunctions in the kernel of GJMS operators. We establish 
that on any manifold of dimension $n\geq 3$, there exist many metrics 
for which our invariants are nontrivial.  We discuss new applications 
to curvature prescription problems.  
\end{abstract}

\maketitle

\section* {Introduction} 
In this note, we announce various results about conformal invariants arising from nodal sets and nullspaces of conformally invariant operators, including some
applications to $Q_{k}$-curvatures. The invariants include nodal sets and the topology of nodal domains
of eigenfunctions in the kernel of GJMS operators (see Section~\ref{section:nodal-sets}). 

We also establish 
that, on any manifold of dimension $n\geq 3$, there exist many metrics 
for which our invariants are nontrivial (see Theorem~\ref{neg:eig:large}). 
In addition, we discuss new applications 
to curvature prescription problems (see Section~\ref{sec:Qk-curvature}). 

This note is organized as follows. In Section~\ref{sec:GJMS}, we
review the main definitions and properties about the GJMS operators
$P_{k,g}$. In Section~\ref{section:nodal-sets}, we describe various
conformally invariant quantities arising from the nullspaces of the GJMS
operators. In Section~\ref{sec:negative-eigv}, we present several
results related to their negative values. In Section~\ref{section:
  heisenberg}, we report on a careful analysis of the spectral properties
of the Yamabe and Paneitz operators on compact Heisenberg
manifolds. In Section~\ref{sec:Qk-curvature}, we study
$Q_{k}$-curvature problems.

Full details on the results mentioned in this note will appear
in~\cite{CGJP1,CGJP2}.

\section{Conformal Powers of the Laplacian and $Q_{k}$-curvatures}\label{sec:GJMS}
Let $(M^{n},g)$ be a Riemannian manifold ($n\geq 3$). We are
interested in the conformal powers of the Laplacian as constructed
in~\cite{GJMS}; these are called the GJMS operators.  

For any positive integer $k$ if $n$ is odd, or for any positive integer~$k\leq
\frac{n}{2}$ if $n$ is even, there is a covariant differential
operator $P_{k,g}:C^{\infty}(M)\rightarrow C^{\infty}(M)$ such that
\begin{enumerate}
    \item[(i)]  $P_{k}$ has with same leading part as $\Delta_{g}^{k}$, where $\Delta_{g}$ 
is the Laplacian on $M$. 

    \item[(ii)] $P_{k}$ is conformally invariant in the sense that
    \begin{equation}
P_{k,e^{2\Upsilon} g}=e^{-(\frac{n}{2}+k)  \Upsilon}P_{k,g}
e^{ (\frac{n}{2}-k) \Upsilon} \qquad \forall \Upsilon\in C^{\infty}(M,\R).
\label{eq:conf-inv-Pk}
\end{equation}
\end{enumerate}

For $k=1$ we get the Yamabe operator $P_{1,g}:=\Delta_{g}+\frac{n-2}{4(n-1)}R_{g}$, where $R_{g}$ is the scalar 
curvature. For $k=2$ we recover the Paneitz operator. In any case, the operator $P_{k}$ is formally self-adjoint (see, 
e.g., \cite{FG02}).

For $w\in \R$ denote by $\cE[w]$ the space of (smooth) conformal
densities of weight $w$ on $M$. A conformal density of weight $w$ can
be defined as a section of the line bundle of conformal densities of
weight $w$ (see, e.g., \cite{PR}). It can be equivalently seen as a
family of $(u_{\hat{g}})_{\hat{g}\in [g]}\subset C^{\infty}(M)$
parametrized by the conformal class $[g]$, in such way that
\begin{equation*}
    u_{e^{2\Upsilon} g}(x)=e^{w\Upsilon(x)}u_{g}(x) \qquad \forall \Upsilon\in C^{\infty}(M,\R).
\end{equation*}Then~(\ref{eq:conf-inv-Pk}) exactly means that $P_{k,g}$ uniquely extends to a differential operator,
\begin{equation}
 P_{k,g}:\cE\biggl[-\frac{n}2+k\biggr]\longrightarrow \cE\biggl[-\frac{n}2-k\biggr],  \label{eq:Pk-densities}
\end{equation}which is independent of the representative metric in the conformal class $[g]$. 

It follows from this that null eigenvectors of $P_{k,g}$ can be
interpreted as conformal densities of weight
$-\frac{n}{2}+k$. Incidentally, the dimension of the nullspace of
$P_{k,g}$ is an invariant of the conformal class $[g]$.

In addition, the operator $P_{k,g}$ is intimately related to the $Q_{k}$-curvature $Q_{k,g}$, which, for $k\neq \frac{n}{2}$, is defined 
by 
\begin{equation*}
    Q_{k,g}:=\frac{2}{n-2k}P_{k,g}(1).
\end{equation*}
For $k=\frac{n}{2}$ the $Q_{\frac{n}{2}}$-curvature is often defined
by a limit procedure by letting $k\rightarrow \frac{n}{2}$
(see~\cite{BO,GZ}). It is often referred to as Branson's
$Q$-curvature.

% % %

\section{Conformal Invariants from the Nullspace of $P_{k,g}$}\label{section:nodal-sets}

Let $k\in \N$ and further assume $k\leq
\frac{n}{2}$ if $n$ is even. If $u=(u_{\hat{g}})_{\hat{g}\in [g]}$ is
a conformal density of weight $w$, then the zero locus
$u^{-1}_{\hat{g}}(0)$ is independent of the metric $\hat{g}$, and
hence is an invariant of the conformal class $[g]$.

In view of~(\ref{eq:Pk-densities}) we 
can regard null-eigenvectors of $P_{k,g}$ as conformal densities of weight $-\frac{n}{2}+k$. Therefore, we obtain

\begin{proposition}\label{prop:nodal-sets-Pk}~     
\begin{enumerate}
        \item If $\dim \ker P_{k,g}\geq 1$, then the nodal set and
          nodal domain of any nonzero null-eigenvector of $P_{k,g}$ give rise to invariants of the conformal class
          $[g]$.
    
        \item If $\dim \ker P_{k,g}\geq 2$, then (non-empty)
          intersections of nodal sets of null-eigenvectors of
          $P_{k,g}$ and their complements are invariants of the
          conformal class~$[g]$.
    \end{enumerate}
\end{proposition}

For a conformal density of weight $0$, all its level sets, not just
its zero-locus, give rise to conformal invariants. Observe that if $n$
is even, then for $k=\frac{n}{2}$ the nullspace of $P_{\frac{n}{2}}$
gives rise to a subspace of $\cE[0]$. It can be shown that the constant functions are contained in the nullspace of 
$P_{\frac{n}{2}}$. We thus obtain

\begin{proposition}
  Assume $n$ is even. If $\dim \ker P_{\frac{n}{2},g}\geq 2$, then the
  level sets of any non-constant null-eigenvector of $P_{\frac{n}{2},g}$ are invariants of the conformal class $[g]$.
\end{proposition}

Let $u_{1,g}, \cdots, u_{m,g}$ be a basis of $\ker P_{k,g}$ (where 
the $u_{j,g}$ are conformal densities of weight $k-\frac{n}{2}$). In addition, set  $\cN:= \bigcap_{1\leq j \leq m} 
u_{j,g}^{-1}(0)$ and define $\Phi:M\setminus \cN \rightarrow \R\proj^{m-1}$ by 
\begin{equation*}
    \Phi(x):=(u_{1,g}(x): \cdots :u_{m,g}(x)) \qquad \forall x \in M\setminus \cN. 
\end{equation*}

Observe that if $x\in M\setminus 
\cN$, then the $m$-uple $(u_{1,g}(x), \cdots u_{m,g}(x))$ depends on the representative metric $g$ only up to positive 
scaling. That is, the projective vector $(u_{1,g}(x): \cdots :u_{m,g}(x))\in \R\proj^{m-1}$ is independent of $g$. 
Therefore, we have

\begin{proposition}\label{prop:embedding-Phi}
   The map $\Phi$ above is an invariant of the conformal class $[g]$. 
\end{proposition}

For $k=\frac{n}{2}$, the nullspace of $P_{\frac{n}{2}}$ always
contains the constant functions, so we may assume that
$u_{1,g}(x)=1$. Moreover, as the $u_{j,g}(x)$ are conformal densities
of weight $0$, for any $x \in M$, the $(m-1)$-uple
$(u_{2,g}(x),\cdots, u_{m,g}(x))$ is actually \emph{independent} of
the representative metric $g$.  Therefore, defining
$\Psi:M\rightarrow \R^{m-1}$ by
\begin{equation*}
   \Psi(x):= (u_{2,g}(x),\cdots ,u_{m,g}(x)) \qquad \forall x \in M,
\end{equation*}we obtain

\begin{proposition}
   The map $\Psi$ above is an invariant of the conformal class $[g]$.   
\end{proposition}

Finally, assume $M$ compact and denote by $dV_{g}(x)$ the Riemannian
measure defined by $g$. If $u_{g}$ is a conformal density of weight
$w<0$ and we set $p= \frac{n}{|w|}$, then a simple exercise shows that
the value of the integral $\int_{M}|u(x)|^{p}dV_{g}(x)$ is independent
of the representative metric $g$. Applying this result to
null-eigenvectors of $P_{k,g}$ we arrive at the following statement.

\begin{proposition}
Assume $M$ compact and $k<\frac{n}{2}$. Let $u_{g}$ be a null-eigenvector of $P_{k,g}$ and regard it as a conformal density of weight 
$-\frac{n}{2}+k$. Then the integral ${\displaystyle
  \int_{M}|u_{g}(x)|^{\frac{2n}{n-2k}}dV_{g}(x)}$ is an invariant of the
conformal class $[g]$.
\end{proposition}

\section{Negative Eigenvalues of the $P_{k,g}$} \label{sec:negative-eigv}
In this section, we assume $M$ compact. Let $k \in \N$ (and further assume 
$k\leq \frac{n}{2}$ if $n$ is even). We are 
interested in metrics for which $P_{k,g}$ has negative eigenvalues. 

For $m \in \N_{0}$, denote by $\cG_{k,m}$ the set of metrics $g$ on
$M$ such that $P_{k,g}$ has at least $m$ negative eigenvalues, where
the eigenvalues are counted with multiplicity. Using the continuity of
the eigenvalues of the $P_{k,g}$ with respect to the metric
(cf.~\cite{Ca}), it is not hard to show that
\begin{itemize}
    \item[-] $\cG_{k,m}$ is an open subset of the space of metrics
      equipped with the $C^{2k}$-topology.

    \item[-] If $g \in \cG_{k,m}$, then the whole conformal
      class~$[g]$ is contained in $\cG_{k,m}$.
\end{itemize}
It follows from this that the number of negative eigenvalues defines a
partition of the set of conformal classes. 

The family $\{\cG_{k,m}\}_{m\geq 0}$ is a non-increasing sequence of open subsets. The question 
that naturally arises 
is whether $\cG_{k,m}$ can be empty for large $m$. 

In the case of the Yamabe operator, using results of
Lokhamp~\cite{Lo96} we get

\begin{theorem}\label{neg:eig:large}  Assume $M$ compact. Then
   \begin{enumerate}
       \item  For any $m\in \N$, there is a metric $g$ on $M$ for which the Yamabe operator 
$P_{1,g}$ has at least $m$ negative eigenvalues. 
     \item There are infinitely many metrics on $M$ for which the
       nullspace of $P_{1,g}$ has dimension~$\geq 1$.
   \end{enumerate}
\end{theorem}

In particular, we see that there are infinitely many metrics that fall into the scope of application of Proposition~\ref{prop:nodal-sets-Pk} and
Proposition~\ref{prop:embedding-Phi}.

It would be interesting to obtain similar results for higher powers of the Laplacian. 
A first result in this direction 
is the following.

\begin{theorem}
Assume $M=\Sigma \times \Sigma$, where $\Sigma$ is a compact surface of genus $\geq 2$. 
Then, for any $m\in \N$, there is
a metric $g$ on $M$ for which the Paneitz operator $P_{2,g}$ has at least $m$ negative 
eigenvalues. 
\end{theorem}
There is a similar result on compact Heisenberg manifolds (see next section).

In addition, as an application of Courant's nodal domain theorem, we obtain 
\begin{theorem}\label{number:nodal}
Let $g$ be a metric such that the Yamabe operator $P_{1,g}$ has exactly $m$ negative 
eigenvalues. Then any (non-zero) 
null-eigenvector has at most $m+1$ nodal domains.  
\end{theorem}

It would be interesting to extend this result to higher conformal powers of the Laplacian.  

\begin{remark}
In contrast to the conformal invariance of $\dim \ker P_{k,g}$, the
non-zero eigenvalues of $P_{k,g}$ exhibit a very different
behavior. As shown by Canzani in~\cite{Ca}, arbitrary small
pertubations within the conformal class can produce simple non-zero
eigenvalues. In fact, if we equip the set of metrics with the
$C^{2k}$-topology, then $P_{k,g}$ has simple non-zero eigenvalues for
all metrics but that of a meager set~(see~\cite{Ca}).
\end{remark}

\section{The Yamabe and Paneitz Operators on Compact Heisenberg Manifolds}\label{section: heisenberg}
In this section, we consider the example of a Heisenberg manifold $M:=\Gamma \backslash \bH_{d}$, $d\in \N$. Here 
$\bH_{d}$ is the $(2d+1)$-dimensional Heisenberg group consisting of unipotent matrices, 
\begin{equation*}
A=    \begin{pmatrix}
1 & x & t\\
0 & 1 & y^{T}\\
0&0&1
\end{pmatrix}, \qquad x,y\in \R^{d}, \quad t\in \R,
\end{equation*}and $\Gamma$ is its lattice subgroup of unipotent matrices with integer-entries. 

In the sequel we shall use coordinates $x=(x_{1},\cdots x_{d})$, $y=(y_{1},\cdots y_{n})$ and $t$ as above to represent an element of
$\bH_{d}$. A left-invariant frame of vector fields on $\bH_{d}$ is then given by the left-invariant vector fields,
\begin{equation*}
 T=\frac{\bd}{\bd t}, \qquad  X_{j}=\frac{\bd}{\bd x_{j}},
\qquad Y_{j}=\frac{\bd}{\bd y_{j}}+x_{j}\frac{\bd}{\bd t},
\quad j=1,\cdots, d.
\end{equation*}As these vector fields are left-invariant, they descend to vector fields on $M$. 

Following~\cite{GW}, for $s>0$, we endow $\bH_{d}$ with the
left-invariant metric,
\begin{equation}
    g_{s}:= \sum_{1\leq j \leq d}dx_{j}
\otimes dx_{j}+\sum_{1\leq j \leq d}s^{-2}dy_{j}\otimes dy_{j}+s^{2d} \theta
    \otimes \theta,
    \label{eq:Heisenberg.metric-gs}
\end{equation}where $\theta:=dt -\sum_{1\leq j \leq d}x_{j}dy_{j}$ is the standard contact form of $\bH_{d}$. This 
metric descends to a metric on $M$. 

For the metric $g_{s}$ the Laplace operator is
$\Delta_{g_{s}}=-\sum_{1\leq j \leq d}(X_{j}^{2}+s^{2}Y_{j}^{2})
-s^{-2d}T^{2}$ and the scalar curvature is $R_{g_{s}}=
-\frac{d}{2}s^{2d+2}$. Therefore, the Yamabe operator is given by
\begin{equation*}
   P_{1,g_{s}}= -\sum_{1\leq j \leq d}(X_{j}^{2}+s^{2}Y_{j}^{2}) -s^{-2d}T^{2}-\frac{2d-1}{16}s^{2d+2}. 
\end{equation*}
The Paneitz operator too can be computed. We find
\begin{multline*}
    P_{2,g_{s}}=\Delta_{g_{s}}^{2}+\frac{12-(2d-1)^{2}}{8(2d-1)}s^{2d+2}\Delta_{g_{s}}+2\frac{d+1}{2d-1}s^{2}T^{2}\\ +
    (2d-3)\frac{(2d+1)(2d-1)^{2}-4(22d+1)}{256(2d-1)^{2}}s^{4d+4}.
\end{multline*}

Using the representation theory of the Heisenberg group enables us to
give explicit spectral resolutions of the Yamabe and Paneitz
operators. In this note, we shall only present some applications of
these spectral resolutions.  A fully detailed account is given
in~\cite{CGJP1}.
\begin{proposition}
    There is a constant $c(d)$ depending only on $d$ such that, for all $s>0$ large enough, the number of negative 
    eigenvalues of the Yamabe operator $P_{1,g_{s}}$ is~$\geq c(d)s^{2d+1}$. 
\end{proposition}
This result provides us with a quantitative version of Theorem~\ref{neg:eig:large} on Heisenberg manifolds. 

In addition, the nullspace of $P_{1,g_{s}}$ can be expressed in terms of Jacobi's theta function,
\begin{equation*}
    \vartheta(z,\tau):= \sum_{k\in \Z} e^{i\pi k^{2}\tau}e^{2i\pi k z},  \qquad z,\tau\in  
    \mathbb{C},\  \Im \tau>0. 
\end{equation*}

\begin{proposition}Suppose that $s^{d+2}=\frac{8\pi}{2d-1}\left(2d+\sqrt{4d^{2}+2d-1}\right)$. Then
    \begin{enumerate}
        \item  The nullspace of the Yamabe operator $P_{1,g_{s}}$ is spanned by the functions,
        \begin{equation*}
            u_{\pm}(x,y,t):=e^{\pm 2i\pi t}e^{\mp \pi s|x|^{2}} \prod_{1\leq j \leq d}\vartheta(y_{j}\pm i sx_{j}, is).
        \end{equation*}
    
        \item  The nodal set of the function $u_{\pm}(x,y,t)$ is equal to the join, 
        \begin{equation*}
           \bigcup_{1\leq j \leq d}\left\{(x,y,t)\in M;  x_{j}=y_{j}=\frac{1}{2}\right\}. 
        \end{equation*}
    \end{enumerate}
\end{proposition}

Finally, for the Paneitz operator in low dimension we obtain
\begin{prop}\label{Paneitz:nil:d=1,2,3} 
For $d=1,2,3$ the number of negative eigenvalues of the Paneitz operator $P_{2,g_{s}}$ goes to $\infty$ as $s\to\infty$.  
\end{prop}

\section{$Q_{k}$-Curvature Prescription Problems} \label{sec:Qk-curvature}
The $Q_{k}$-curvature prescription problem has been an
important focus of interest in conformal geometry, especially for Branson's $Q$-curvature $Q_{\frac{n}{2}}$ ($n$ even) 
(see, e.g., \cite{BaFaRe,Brendle,DM}).  
In this section, we look at the
$Q_{k}$-curvature prescription problem for $k\neq \frac{n}{2}$
(further assuming $k<\frac{n}{2}$ when $n$ is even).

In the sequel, we denote by $\cR(Q_{k})$ the set of functions $v \in C^{\infty}(M)$ for which there is a metric $\hat{g}$ in the 
conformal class $[g]$ such that $v=Q_{k,\hat{g}}$. 

Let $\hat{g}$ be a metric in the conformal class $[g]$. The definition of $Q_{k,\hat{g}}$ implies that 
$Q_{k,\hat{g}}=0$ if and only if $P_{k,\hat{g}}(1)=0$. Furthermore, it follows from the 
transformation law~(\ref{eq:conf-inv-Pk}) that $P_{k,\hat{g}}(1)=0$ if and only if the nullspace of $P_{k,g}$ contains a positive 
function. Therefore, we obtain

\begin{theorem}\label{Pk:kernel:const}
$\cR(Q_{k})$ contains the zero function if and only if the nullspace
  of $P_{k,g}$ contains a nowhere-vanishing function.
\end{theorem}

For the Yamabe operator $P_{1,g}$ we actually get a finer result. 
\begin{theorem}\label{int:negative}
Let $u_{g}\in \ker P_{1,g}\setminus \{0\}$ and let $\Omega$ be a nodal domain of $u_{g}$. Then, for any metric $\hat{g}$ in the 
conformal class $[g]$, 
\begin{enumerate}
    \item ${\displaystyle \int_\Omega |u_{\hat{g}}(x)|\, R_{\hat g}(x)dV_{\hat g}(x) <0}$, where $R_{\hat{g}}$ is the 
    scalar curvature.

    \item  The scalar curvature $R_{\hat{g}}$ cannot be nonnegative on $\Omega$. 
\end{enumerate}     
\end{theorem}

Finally, assume $M$ is compact. In \cite{Gov10} constraints on 
$\mathcal{R}(Q_{\frac{n}{2}})$ are described; these arise from 
$\ker P_{\frac{n}{2}}$.  In Appendix to~\cite{CGJP1},
Gover-Malchiodi prove the following.

\begin{theorem}\label{prop:vanishing}
    Let $v \in \cR(Q_{k})$. Then, for all $u\in
    \ker P_{k,g}\setminus \{0\}$, there is a metric $\hat{g}$ in the
    conformal class $[g]$ such that
    \begin{equation*}
        \int_{M}u(x)v(x)dV_{\hat{g}}(x)=0.
    \end{equation*}
\end{theorem}

When $\ker P_{k,g}$ is non-trivial this results yields infinitely many
constraints on $\cR(Q_{k})$. In particular, there is an
infinite-dimensional (non-linear) subspace of $C^{\infty}(M)$ that is
disjoint from $\cR(Q_{k})$.
  
An immediate consequence of Theorem~\ref{prop:vanishing} is the following. 

\begin{theorem}
  Let $v \in C^{\infty}(M)$ and assume there is $u\in \ker
  P_{k,g}\setminus \{0\}$ such that $uv$ is non-negative everywhere
  and positive at a point. Then $v$ cannot be contained in
  $\cR(Q_{k})$.
\end{theorem}

Gover-Malchiodi apply the above result to two special cases:
\begin{itemize}
    \item[-] $u\in \ker P_{k,g}\setminus \{0\}$ and $v=s_u$ where
      $s_u$ is any function with the same strict sign as $u$.

    \item[-]  $u\in \ker P_{k,g}\setminus \{0\}$ is non-negative (resp., non-positive) and $v$ is a positive (resp., 
    negative) constant function.  
\end{itemize}
As a result, they obtain the following:

\begin{theorem} Assume $\ker P_{k,g}$ is non-trivial.   
\begin{enumerate}
    \item  
For $u \in \ker P_{k,g}\setminus \{0\}$ no function with the same
strict sign as $u$ can be contained in $\cR(Q_{k})$. In particular
$u\notin \cR(Q_{k})$.

    \item  If $\ker P_{k,g}\setminus \{0\}$ contains a non-negative function, then there is no metric 
     in the conformal class $[g]$ with non-zero constant $Q_{k}$-curvature.
\end{enumerate}
\end{theorem}

\end{document}